\documentclass[11pt]{article}
\usepackage{}
\usepackage{amssymb,amsmath,verbatim,amsthm}
\newtheorem{Theorem}{Theorem}[section]
\newtheorem{Lemma}[Theorem]{Lemma}
\newtheorem{Corollary}[Theorem]{Corollary}

\newtheorem{Example}[Theorem]{Example}

\newtheorem{Proposition}[Theorem]{Proposition}
\newtheoremstyle{noparens}
  {}{}
  {\itshape}{}
  {\bfseries}{.}
  { }
  {\thmname{#1}\thmnumber{ #2}\mdseries\thmnote{ #3}}
\theoremstyle{noparens}
\usepackage{enumerate}
\usepackage{color,xcolor}
\usepackage{indentfirst}
\newcommand{\ignore}[1]{}

\begin{document}

\renewcommand{\theequation}{\thesection.\arabic{equation}}

\title{\bf Tilings of $\mathcal{H}_{q}(n,w)$ with optimal $(n,d,w)_{q}$-codes\footnote{Supported
 by National Natural Science Foundation of China (12571346, 12371326).}}
 \author{
 { Yuli  Tan,}  { Junling  Zhou\footnote{Corresponding Author}}\\
 {\small School of Mathematics and Statistics}\\ {\small Beijing Jiaotong University}\\
  {\small Beijing  100044, China}\\
 {\small YL.Tan@bjtu.edu.cn}\\
{\small jlzhou@bjtu.edu.cn}\\
}
\date{ }
\maketitle

\begin{abstract}
The metric space $\mathcal{H}_{q}(n,w)$ is the set of all words of length $n$ with weight $w$ over the alphabet $\mathbb{Z}_{q}$, under the Hamming distance metric.
A $q$-ary constant-weight code, as a nonempty subset of $\mathcal{H}_{q}(n,w)$, has always been a fundamental topic in coding theory. This paper investigates the tiling problem of $\mathcal{H}_{q}(n,w)$ with optimal $(n,d,w)_{q}$-codes, simply denoted by $\mathrm{TOC}_{q}(n,d,w)$, meaning a partition of $\mathcal{H}_{q}(n,w)$ into mutually disjoint optimal $q$-ary constant-weight codes  with distance $d$. When the distance $d$ is odd, we investigate large sets of generalized Steiner systems. When $d$ is even, we define large sets of generalized maximum H-packings. We present several general construction approaches for generating $\mathrm{TOC}_{q}(n,d,w)$s via  $t$-resolvable Steiner systems and almost-regular edge-colorings of complete hypergraphs. For the cases $d=2$ and $d=2w$, we completely resolve the existence problem of $\mathrm{TOC}_{q}(n,d,w)$s for all parameters $q,n$ and  $w$. Particularly, we pay attention to tilings for weight three. For binary case and weight three, the existence problem of  $\mathrm{TOC}_{2}(n,d,3)$s is  totally resolved.  For specific alphabet size $q\ge 3$, we obtain many infinite families of $\mathrm{TOC}_{q}(n,d,3)$s for distances $d=3,4,5$.

\medskip

\medskip
\end{abstract}

\section{Introduction}

Constant-weight codes (CWCs) are an important and fascinating class of codes. CWCs with a fixed distance possess error detection and error correction capabilities \cite{code-impo1}.
Binary CWCs are a canonical research subject in coding theory, due to their applications in asymmetric channels \cite{[16], [39], [51], [90]}, DC-free codes \cite{[64], [175]},  spherical codes \cite{[167]}, and recently in storage systems \cite{[6]}.
During  the past three decades,  $q$-ary CWCs with  $q>2$ aroused much interest among researchers in coding theory and combinatorics, primarily because of their increasing applications in powerline communications \cite{[2]}, DNA computing  \cite{apply2,apply3}, coding for bandwidth-efficient channel \cite{apply1}.

The fundamental problem in coding theory is that of determining the value $A_{q}(n,d,w)$, that is, the maximum size of a $q$-ary code of length $n$ with constant weight $w$ and distance $d$. When $q=2$, the subscript $q$ is usually omitted. Obviously the Hamming distance between two codewords in a binary CWC is even. The nature of
$A(n,d,w)$ is combinatorial, because this value equals the so-called packing number in design theory, namely,  $A(n,2w - 2t + 2,w) =D(t,w,n)$, where the packing number $D(t,w,n)$ denotes the
maximum number of $w$-subsets (called blocks) of a $n$-set $S$, such that every
$t$-subset of $S$ is contained in at most one block. By this equivalence, for binary CWCs, $A(n,d,w)$ has been widely studied;  the cases of $w=3,4$ have been totally determined \cite{Bao-Ji,CRC} (certainly $w=2$ is trivial).

For alphabet size $q>2$,  optimal $q$-ary CWCs (containing $A_{q}(n,d,w)$ codewords) are also closely related with combinatorial objects.  Johnson-type bounds  on  $A_{q}(n,d,w)$ were established in \cite{bound_d2, Svanstrom1999-1}. Then
 a large literature suggests that many combinatorial configurations play effective roles on constructions of optimal $q$-ary CWCs, such as generalized Steiner systems \cite{TUVI1}, group divisible designs with restrictions \cite{Ge2002-1, Ge2002-2}, completely reducible super-simple designs \cite{wu, Wei}, large sets of designs \cite{Chee-2008}. For the cases $d=2,2w$ or $w=2$, the determination of $A_{q}(n,d,w)$ is not very difficult and  has been accomplished \cite{Chee-2007, bound_d2}.
 For weight $w=3$,
the exact value of $A_{q}(n,4,3)$ was successfully  determined by Chee et al. \cite{Chee-2008, Chee-2007}.
Subsequently in 2015, Chee et al. \cite{Chee-2015} gave a complete solution to the determination of $A_{q}(n,5,3)$.
For distance $d=3$, there is a large amount of work  focusing on $A_{q}(n,3,3)$,
see for example \cite{Wilson1999, Chen1999, TUVI1, Ge2000, Ge2002-1, Ge2002-2,
Ge2003-1, Phelps1997, Phelps1999, Wu2000, Yin1999}. This problem has  been completely resolved for small alphabet sizes (say, $q\le 13$), although we do not have so satisfactory solutions as for the distances $d=4,5$.
For weight $w=4$, the investigation on $q$-ary CWCs is not quite extensive, see Cao et al. \cite{Cao-2007} for $d=3$, Ge and Wu \cite{Ge2003-2} and Chen et al. \cite{wu} for $d=4$, and Ji et al. \cite{Ji-2005-dcc} for $d=5$.
Recently, Liu and Shangguan \cite{Shangguan1} showed that for all fixed odd distances, there exist near-optimal constant-weight codes asymptotically achieving the classic
 Johnson-type upper bounds. As a complement, Bennett \cite{Bennett2025} resolved the situation where the distance is even.

As we know, a  $q$-ary code of length $n$ with constant weight $w$ is a nonempty subset of the space $\mathcal{H}_{q}(n,w)$, equipped with the Hamming metric, consisting of all words of length $n$
with weight $w$ over a $q$-element alphabet. As previously stated, there has been abundant research  on constructions of optimal
$q$-ary constant-weight codes with given parameters, especially for weight $w=3$. A question arises naturally: Does there exist a tiling  of the  entire space
$\mathcal{H}_{q}(n,w)$ with mutually disjoint optimal $q$-ary CWCs with a fixed distance?  In other words, is there a partition of $\mathcal{H}_{q}(n,w)$ into
 optimal $q$-ary CWCs with a fixed distance? This problem has received little attention. However, we have many analogs under the name of ``large sets" in
 design theory, for example, the earliest large set problem for Kirkman triple systems, proposed by Sylvester in 1850's. In this paper, we will discuss such a large set problem for optimal $q$-ary CWCs. To be more precise, we treat
 the existence problem of a tiling $\mathrm{TOC}_{q}(n,d,w)$, with TOC denoting Tiling of Optimal Codes and $q,n,d,w$ the parameters of the CWC. As the main results, considering all possible parameters,
we first propose several construction approaches of  tilings for general weights  and then display many results of tilings for specific weight three.

The outline of the paper is as follows. In Section 2, we introduce the definitions of constant-weight codes and tilings. The tools in combinational designs and hypergraphs used in this paper are also explored.
In Section 3, we present some general constructions to generate $\mathrm{TOC}_{q}(n,d,w)$s.
In Section 4, we focus on weight three and produce many existence results of $\mathrm{TOC}_{q}(n,d,3)$s with $d=3, 4, 5$.
In Section 5, a brief conclusion and some open problems will be given.

\section{Preliminaries}
In this section, we introduce the definitions of constant-weight codes and tilings. We also recall some basic notations related to combinational designs and hypergraphs.

\subsection{Constant-weight codes and tilings}

Let $a,b\in \mathbb{Z}^{+}$ such that $a<b$. The closed integer interval $[a,b]$ is defined as:
\begin{equation*}\label{E00}
[a,b]=\{x\in\mathbb{Z}^{+}:a\leq x \leq b\}.
\end{equation*}
Specifically, the interval $[1,n]$ is abbreviated as $[n]$ for $n\in \mathbb{Z}^{+}$.
For a positive integer $q$, denote the ring $\mathbb{Z}/q\mathbb{Z}$ by $\mathbb{Z}_{q}$.
Let $\mathbb{Z}_{q}^{n}$ be the set of all words (or vectors) of length $n$ over the alphabet $\mathbb{Z}_{q}$.
A word $\mathbf{x}\in \mathbb{Z}_{q}^{n}$ is denoted by $\mathbf{x} = (x_{1},x_{2},\ldots,x_{n})$.
Let $\mathcal{H}_{q}(n)=\mathbb{Z}_{q}^{n}$ denote the {\em$q$-ary Hamming $n$-space} which is endowed with {\em Hamming distance metric} $d_{H}$ defined as follows:
$$d_{H}(\mathbf{x},\mathbf{y})= |\{i\in [n]:x_{i}\neq y_{i}\}|,\ \mathbf{x},\mathbf{y}\in \mathcal{H}_{q}(n),$$
i.e., the number of coordinates where $\mathbf{x}$ and $\mathbf{y}$ differ. A {\em $q$-ary code} of length $n$ is a subset $\mathcal{C}$ of $\mathcal{H}_{q}(n)$. Every element of $\mathcal{C}$ is called a {\em{codeword}}.
The {\em (Hamming) weight} of a word $\mathbf{x}$ is the quantity $d_{H}(\mathbf{x},\mathbf{0})$.
The {\em support} of $\mathbf{x}$, denoted by $\text{supp}(\mathbf{x})$, is the set of coordinate positions $i\in[n]$ whose entry $x_{i}$ is nonzero. In other words, the weight of a word is the size of its support.

The metric space $\mathcal{H}_{q}(n,w)$ is the set of all words of weight $w$ in $\mathcal{H}_{q}(n)$, under the Hamming distance metric.
A {\em $q$-ary constant-weight code} $\mathcal{C}$ of length $n$, weight $w$ and distance $d$, denoted $(n,d,w)_{q}$-code, is a nonempty subset of $\mathcal{H}_{q}(n,w)$ such that $d_{H}(\mathbf{x},\mathbf{y})\geq d$ for all distinct $\mathbf{x},\mathbf{y}\in \mathcal{C}$.
The number of codewords in an $(n,d,w)_{q}$-code is called the {\em size} of the code. The maximum size of an $(n,d,w)_{q}$-code is denoted $A_{q}(n,d,w)$. If the subscript $q$ is omitted, we mean $q=2$. An {\em optimal} $(n,d,w)_{q}$-code is an $(n,d,w)_{q}$-code having $A_{q}(n,d,w)$ codewords.
Below are classic results regarding the upper bound on the quantity $A_{q}(n,d,w)$.

\begin{Lemma}[\label{classic-bound}Svanstr\"{o}m \cite{Svanstrom1999-1}]
$$A_{q}(n,d,w)\leq \lfloor\frac{n}{n-w}A_{q}(n-1,d,w)\rfloor,$$
$$A_{q}(n,d,w)\leq \lfloor\frac{n(q-1)}{w}A_{q}(n-1,d,w-1)\rfloor.$$
\end{Lemma}

Applying Lemma \ref{classic-bound} recursively yields the following bound.
\begin{Lemma}[\label{classic-bound-generalized}Fu et al. \cite{bound_d2}]
Let $t=\lceil\frac{2w-d+1}{2}\rceil$. Then
\begin{equation}\label{E-2}
A_{q}(n,d,w)\leq
\begin{cases}
\frac{{n\choose t}(q-1)^{t}}{{w\choose t}},&\text{ if }d\text{ is odd}, \\
\frac{{n\choose t}(q-1)^{t-1}}{{w\choose t}},& \text{ if }d\text{ is even}.\\
\end{cases}
\end{equation}
\end{Lemma}

A {\em{tiling}} of the space $\mathcal{H}_{q}(n,w)$ with optimal constant-weight codes of length $n$, weight $w$ and distance $d$, denoted by $\mathrm{TOC}_{q}(n,d,w)$, is a collection of mutually disjoint optimal $(n,d,w)_{q}$-codes such that their union covers the whole space $\mathcal{H}_{q}(n,w)$. Every optimal $(n,d,w)_{q}$-code, as a member of the tiling, is called a {\em{tile}}. We use $m_{q}(n,d,w)$ to represent the number of tiles contained in a $\mathrm{TOC}_{q}(n,d,w)$.
The necessary condition for the existence of a $\mathrm{TOC}_{q}(n,d,w)$ is $A_{q}(n,d,w)\mid{{n}\choose{w}}(q-1)^{w}$. To avoid trivialities, we always assume that $d\geq 2$. In this paper, we place our focus on $3\leq w \leq n$.
According to Lemma \ref{classic-bound-generalized}, we obtain the following lemma.
\begin{Lemma}\label{mq1}
Let $t=\lceil\frac{2w-d+1}{2}\rceil$. If (\ref{E-2}) holds with equality, then
$$m_{q}(n,d,w)=
\begin{cases}
{{n-t}\choose {w-t}}(q-1)^{w-t},&\text{ if }d\text{ is odd}, \\
{{n-t}\choose {w-t}}(q-1)^{w-t+1},& \text{ if }d\text{ is even}.\\
\end{cases}
$$
\end{Lemma}

Now, we give an example of tilings. Note that the upper bound in (\ref{classic-bound}) is not always tight.

\begin{Example}\label{(4,4,3)_3}
{\rm We construct a $\mathrm{TOC}_{3}(4,4,3)$, which consists of $16$ mutually disjoint optimal $(4,4,3)_3$ codes of size $2$.
For $i,j,k\in [2]$, let
$$\mathcal{C}_{i,j,k}^1= \big\{(i\ j\ k\ 0),(0\ \bar{j}\ \bar{k}\ \bar{i})\big\},$$
$$\mathcal{C}_{i,j,k}^2= \big\{(i\ j\ 0\ k),(\bar{i}\ 0\ \bar{j}\ \bar{k}) \big\},$$
where $\bar{1}=2,\bar{2}=1$. Then the collection  $\big\{\mathcal{C}_{i,j,k}^1,\mathcal{C}_{i,j,k}^2:i,j,k\in [2]\big\}$ is a $\mathrm{TOC}_{3}(4,4,3)$.
}
\end{Example}

For convenience, we also represent a word $\mathbf{x}=(x_{1},x_{2},\ldots,x_{n})\in\mathbb{Z}_{q}^{n}$ using a set notation $\{(i,x_{i}):i\in[n],x_{i}\neq 0\}$.
Hence, a word $\mathbf{x}=(x_{1},x_{2},\ldots,x_{n})\in \mathcal{H}_{q}(n,w)$ is also regarded as a set $\{(i_{1},x_{1}),(i_{2},x_{2}),\ldots,(i_{w},x_{w})\}\subseteq[n]\times[q-1]$, where $\{i_{1},i_{2},\ldots,i_{w}\}$ is its support.
When a word has relatively large length but relatively small weight, its set representation is more concise and simpler. By abuse of notation, we do not distinguish between the vector  representation and the set  representation.

\subsection{H-designs, resolvable Steiner systems and orthogonal arrays}\label{H+LS}

Let $X$ be an $ng$-set and $\mathcal{G}$ be a partition of $X$ into $n$ sets of cardinality $g$, called {\em{groups or holes}}.
An {\em $\mathrm{H}$-design $\mathrm{H}(n,g,w,t)$} is a triple $(X,\mathcal{G},\mathcal{A})$ with the following properties:\vspace{-0.2cm}
\begin{itemize}
  \item [(1)] $\mathcal{A}$ is a family of $w$-subsets of $X$, called {\em{blocks}};\vspace{-0.2cm}
 \item [(2)] each block of $\mathcal{A}$ intersects any given group in at most one point;\vspace{-0.2cm}
 \item [(3)] each $t$-subset of $X$ from $t$ distinct groups is contained in exactly one block.
\end{itemize}
\vspace{-0.2cm}
If the group size $g=1$ and $(X,\mathcal{G},\mathcal{A})$ is an $\mathrm{H}(n,1,w,t)$ design, then we omit the notation of $\mathcal{G}$ and call $(X,\mathcal{A})$ a {\em{Steiner system}}, also denoted $\mathrm{S}(t,w,n)$.
An $\mathrm{S}(2,3,n)$ is a {\em{Steiner triple system}} $\mathrm{STS}(n)$ and an $\mathrm{S}(3,4,n)$ is a {\em{Steiner quadruple system}} $\mathrm{SQS}(n)$.

The notion of H-designs has a natural extension to H-packings.
Let the set $X$ and the group set $\mathcal{G}$ be the same as above. An {\em {$\mathrm{H}$-packing}} $\mathrm{HP}(n,g,w,t)$ is a triple $(X,\mathcal{G},\mathcal{A})$, where the block set $\mathcal{A}$ possesses the properties (1), (2) as above and also the property:\vspace{-0.2cm}
\begin{itemize}
 \item [(3')] each $t$-subset of $X$ from $t$ distinct groups is contained in at most one block.
\end{itemize}
\vspace{-0.2cm}
If $(X,\mathcal{G},\mathcal{A})$ is an $\mathrm{HP}(n,1,w,t)$, then we call $(X,\mathcal{A})$ a {\em{$t$-packing}}, denoted $\mathrm{P}(t,w,n)$.
A $t$-packing  $(X,\mathcal{A})$ is called {\em optimal} or {\em maximum} if there does not exist any $t$-packing  $(X,\mathcal{B})$ (with the same parameters) such that $|\mathcal{A}|<|\mathcal{B}|$. Obviously, a Steiner system  $\mathrm{S}(t,w,n)$ is an optimal $t$-packing  $\mathrm{P}(t,w,n)$.

Suppose that $(X,\mathcal{G},\mathcal{A})$ is an $\mathrm{HP}(n,g,w,t)$ where $X=[n]\times[g]$, $\mathcal{G}=\{G_{i}:1\leq i \leq n\}$ with $G_{i}=\{i\}\times[g]$.
Then we may produce a CWC over $\mathbb{Z}_q$ with $q=g+1$.
Any block $A=\{(i_{1},j_{1}),(i_{2},j_{2}),\ldots,(i_{w},j_{w})\}\in \mathcal{A}$ is a set representation of a codeword $C_A$, whose $i_{s}$-th coordinate is $j_{s}$ with $1\leq s \leq w$, and all other coordinates of $C_A$ are zeros.
Then $\mathcal{C}=\{C_A:A\in \mathcal{A}\}$ forms a $q$-ary constant-weight code, called the code of the design $(X,\mathcal{G},\mathcal{A})$.
First we consider the case $g=1$ and thus $(X,\mathcal{A})$ is a  packing $\mathrm{P}(t,w,n)$. Now we obtain a binary constant-weight code from  $(X,\mathcal{A})$.
Since the intersection of two distinct blocks in $\mathcal{A}$ has at most $t-1$ points, it follows that
the distance $d$ of ${\cal C}$ is $2(w-t+1)$. Further if $(X,\mathcal{A})$ is  an optimal packing $\mathrm{P}(t,w,n)$,
then the binary constant-weight code from $(X,\mathcal{A})$ becomes an optimal $(n,2(w-t+1),w)_{2}$-code.
On the other hand, letting $\mathcal{C}$ be an optimal $(n,2(w-t+1),w)_{2}$-code,
the  set of supports $\{\text{supp}(\mathbf{x}):\mathbf{x}\in \mathcal{C}\}$ is the block set of an optimal packing $\mathrm{P}(t,w,n)$. Therefore, an optimal $(n,2(w-t+1),w)_{2}$-code is equivalent to an  optimal $\mathrm{P}(t,w,n)$. This is the well-known equivalence during the process of investigation of binary CWCs.
Now we turn to the case $g\geq 2$.  Whenever $(X,\mathcal{G},\mathcal{A})$ is an $\mathrm{H}$-design $\mathrm{H}(n,g,w,t)$, the associated code $\mathcal{C}$ has  distance $d\le2(w-t)+1$. If the equality holds, i.e., $d=2(w-t)+1$, then we call such an $\mathrm{H}$-design a {\em generalized Steiner system}, denoted $\mathrm{GS}(t,w,n,g)$.
This suggests that  a  $\mathrm{GS}(t,w,n,g)$ gives rise to an optimal $(n,2(w-t)+1,w)_{q}$-code. The converse also holds, as stated in the following lemma.

\begin{Lemma}[\label{Boundgs}{\cite[Corollary 6.1]{Chee-2010}}]
The existence of a $\mathrm{GS}(t,w,n,q-1)$ is equivalent to the existence of an $(n,2(w-t)+1,w)_{q}$-code with $\frac{{n\choose t}(q-1)^{t}}{{w\choose t}}$ codewords.
\end{Lemma}

 The equivalence in Lemma \ref{Boundgs} has been frequently utilized to handle nonbinary CWCs with odd distances.   While for nonbinary CWCs with even distances, $\mathrm{H}$-packings will take effect.
When $q\ge 3$ and the distance $d$ is even, let $t=\frac{2w-d+2}{2}$ and
 assume that there exists an $(n,d,w)_q$-code ${\cal C}$ with  distance $d$. Making use of the set representations of codewords in ${\cal C}$, we arrive at an
 $\mathrm{HP}(n,g,w,t)$. Indeed, it is an  $\mathrm{H}$-packing instead of an  $\mathrm{H}$-design; otherwise, the minimum distance would be at most $2(w-t)+1$,
 as we explored above. 
  We define a {\em generalized maximum $\mathrm{H}$-packing}, denoted $\mathrm{GMHP}(t,w,n,g)$, to be the largest $\mathrm{HP}(n,g,w,t)$ with minimum distance $2(w-t+1)$. Analogous to the  odd distance case, an optimal  $(n,d,w)_q$-code with  even distance $d=2(w-t+1)$ is the same as a  $\mathrm{GMHP}(t,w,n,g)$.

Two $\mathrm{S}(t,w,n)$s on the same point set are said to be {\em{disjoint}} if they have no blocks in common.
An $\mathrm{S}(t, w, n)$ $(X,\mathcal{A})$ is called {\em $i$-resolvable}, $0 <i\leq t$, if its block set $\mathcal{A}$ can be partitioned into pairwise disjoint $\mathrm{S}(i, w, n)$s. Such a partition is called an $i$-resolution.
A $1$-resolvable $\mathrm{S}(t, w, n)$ is called {\em resolvable} and denoted $\mathrm{RS}(t, w, n)$, where each $\mathrm{S}(1, w, n)$ in the $1$-resolution is called a {\em parallel class} of $X$.
An $\mathrm{RS}(2,3,n)$ is a {\em{Kirkman triple system}} $\mathrm{KTS}(n)$.
It is well-known that a $\mathrm{KTS}(n)$ exists if and only if $n\equiv  3 \pmod{6}$. 

\begin{Theorem}\label{KTS}
\textup{\rm(\cite{3-RS})} Let $q$ be a prime power and $n$ be a positive integer. Then there exists an $\mathrm{RS}(3,q + 1,q^{n} + 1)$ if and only if $n$ is odd.
\end{Theorem}

A $2$-resolvable $\mathrm{SQS}(n)$ is comprised of $\frac{n-2}{2}$ mutually disjoint $\mathrm{S}(2,4,n)$s. Hence, the necessary condition for the existence of a $2$-resolvable $\mathrm{SQS}(n)$
is $n\equiv 4\pmod {12}$. Now, we recall some known results.

\begin{Theorem}[\label{2RSQS1}{\cite{Baker-1976,Teirlinck-1994}}]
There exists a $2$-resolvable $\mathrm{SQS}(n)$ if $n= 4^{m}$ or $n=2p^{m}+2$ with $p\in\{7,31,127\}$ for any positive integer $m$.
\end{Theorem}

A $t$-resolvable $\mathrm{S}(w,w,n)$ is well-known as a {\em{large set of Steiner systems}} and denoted by $\mathrm{LS}(t,w,n)$. So an $\mathrm{LS}(t,w,n)$ is a partition of all $w$-subsets of a $n$-set $X$ into pairwise disjoint block sets $\mathcal{B}_{i}$ ($1\leq i \leq m$) such that each $(X,\mathcal{B}_{i})$ is an $\mathrm{S}(t,w,n)$. Obviously, $m={{n-t}\choose {w-t}}$.
It is well-known that a large set of $\mathrm{STS}(n)$s exists if and only if
$n \equiv 1,3\pmod{6}$, $n\geq 3$ and $n\neq 7$ \cite{Lu13,Lu46,Teirlinck L}.
In an $\mathrm{LS}(t,w,n)$, if each member $\mathrm{S}(t,w,n)$ is resolvable, we call it a {\em{large set of resolvable Steiner systems}}, denoted $\mathrm{LRS}(t,w,n)$.
An $\mathrm{LRS}(2,3,n)$ is a {\em{large set of Kirkman triple systems}}, denoted $\mathrm{LKTS}(v)$.
The problem of the existence of $\mathrm{LKTS}$s is an important problem in combinatorial design theory since 1850.
Although it has attracted considerable attention from researchers and has been the subject of numerous studies, it has been far from completely resolved and remains wide open up to now.
A tripling construction states that an $\mathrm{LKTS}(n)$ yields an $\mathrm{LKTS}(3n)$ \cite{zhang-lkts}.

We introduce the notion of {\em{overlarge set of Kirkman triple systems}}, a similar design to $\mathrm{LKTS}$, denoted $\mathrm{OLKTS}(n)$.
It is a partition of all triples of an $(n+1)$-set $X$ into pairwise disjoint block sets $\mathcal{B}_{x}$ ($x\in X$) such that each $(X\setminus\{x\},\mathcal{B}_{x})$
is a $\mathrm{KTS}(n)$. 
There are many works studying the existence of $\mathrm{OLKTS}$s \cite{lkts+olkts,olkts1,olkts4,olkts2,olkts3}.
An $\mathrm{RDSQS}(n+1)$ is a Steiner quadruple system $\mathrm{SQS}(n+1)$ whose derived designs are resolvable at all points, meaning that an  $\mathrm{RDSQS}(n+1)$ gives an  $\mathrm{OLKTS}(n)$. A tripling construction states that an $\mathrm{RDSQS}(n+1)$ yields an $\mathrm{RDSQS}(3n+1)$ \cite{yuan-rdsqs}.
In this paper, we omit a comprehensive list of the results about $\mathrm{LKTS}$s and $\mathrm{OLKTS}$s, but we display some results derived from combining the tripling constructions for LKTSs and RDSQSs with known small orders less than $100$.

\begin{Theorem}[\label{Lkts}\cite{lkts+olkts,C-rd,LEI1,Tan1,lkts-yuan,yuan-rdsqs,zhou1,zheng,zhang-lkts}]
Let $m$ be any nonnegative integer. For any positive integer $a\equiv 3\pmod{6}$ and $a<100$, there exist
\begin{itemize}
  \item [$(1)$]  an $\mathrm{LKTS}(3^{m}\cdot a)$ if $a\not\in \{57,87,93\}$ and \vspace{-0.2cm}
  \item [$(2)$]  an $\mathrm{OLKTS}(3^{m}\cdot a)$ if $a\not\in \{51,69,75\}$. \vspace{-0.2cm}
  \end{itemize}
\end{Theorem}


An {\em orthogonal array} $\mathrm{OA}(t,n,q)$ is a $q^{t}\times n$ array $A$ over an alphabet $Q$ with $q$ symbols such that each ordered $t$-tuple of $Q$ appears in exactly one row in each projection of $t$ columns from $A$. 
The following results \cite[pp. 317--331]{code-impo1} are derived from $\mathrm{MDS}$ codes.
\begin{Theorem}[\label{OA_B4}{\cite{code-impo1}}]
If $q$ is a prime power and $t$ is a positive integer with $t\leq q-1$, then there exists an $\mathrm{OA}(t,n,q)$ for any integer $n$ with $t\leq n\leq q+1$. Moreover, if $q$ is an even prime power and $t\in \{3,q-1\}$, then there exists an $\mathrm{OA}(t,n,q)$ for any $n$ with $t\leq n\leq q+2$.
\end{Theorem}
%
%


\subsection{Hypergraphs, almost-regular edge-colorings and strong vertex-colorings}

A {\em hypergraph} $\mathcal{H}$ is a pair $(V(\mathcal{H}),E(\mathcal{H}))$, where $V(\mathcal{H})$ is a set of {\em vertices} and $E(\mathcal{H})$ is a collection of subsets of $V(\mathcal{H})$, called {\em edges}.
The {\em{rank}} of $\mathcal{H}$ is $r(\mathcal{H})=\mathop{\max}\{|E|:E\in E(\mathcal{H})\}$. If $r(\mathcal{H})=2$, $\mathcal{H}$ is called a {\em{graph}}.
If all edges have the same size $w$, we say that $\mathcal{H}$ is a {\em{$w$-uniform}} hypergraph.
For any positive integer $n$, the {\em{complete $w$-uniform hypergraph}} $K_{n}^{w}$ is the hypergraph  comprising of all $w$-subsets of  $n$ vertices as edges. 

Let $\mathcal{H}$ be a hypergraph. The {\em{degree}} of $x\in V(\mathcal{H})$ in $\mathcal{H}$ is $\text{deg}(x)=|\{x\in E:E\in E(\mathcal{H})\}|$, i.e., the number of edges that contain $x$ in the edge set $E(\mathcal{H})$.
If $\text{deg}(x)= r$ for all $x\in V(\mathcal{H})$, we say that $\mathcal{H}$ is {\em{$r$-regular}}.
A {\em{$1$-factor}} of $\mathcal{H}$ is a $1$-regular spanning sub-hypergraph of $\mathcal{H}$, and a {\em{near $1$-factor}} of $\mathcal{H}$ is a sub-hypergraph having one isolated vertex of degree zero and  all other vertices degree $1$.
A {\em{$1$-factorization}} or {\em{near $1$-factorization}} of $\mathcal{H}$ is a set of $1$-factors or near $1$-factors, that partitions the edges of $\mathcal{H}$, respectively.
A hypergraph is {\em{almost-regular}} if
$|\text{deg}(x) - \text{deg}(y)|\leq 1$ for any two vertices $x$ and $y$.
An {\em{almost-regular edge-coloring}} of $\mathcal{H}$ is a partition of the edges in $E(\mathcal{H})$ into color classes, where each color class forms an almost-regular sub-hypergraph of $\mathcal{H}$. Hence, $1$-factorization and near $1$-factorization are all almost-regular edge-colorings.
The following theorem is from Baranyai \cite{Baranyai}.

\begin{Theorem}[\label{almostregular}\cite{Baranyai}]
Let $n,w$ be positive integers with $\lfloor\frac{n}{w}\rfloor\mid{{n}\choose{w}}$.
Then the complete $w$-uniform hypergraph $K_{n}^{w}$ has an almost-regular edge-coloring with ${{n}\choose{w}}/{\lfloor\frac{n}{w}\rfloor}$ color classes such that the number of edges in each color class is $\lfloor\frac{n}{w}\rfloor$.
\end{Theorem}

As a corollary to Theorem \ref{almostregular}, a $1$-factorization of $K_{n}^{w}$ exists if and only if $w\mid n$ and a near $1$-factorization of $K_{n}^{w}$ exists if and only if $n\equiv 1\pmod{w}$ and $\frac{n-1}{w}\mid{n\choose w}$.

The hypergraph $\lambda K_{n}^{w}$ has vertex set $[n]$ and edges given by including $\lambda$ copies of each $w$-subset of $[n]$. When $\lambda =1$, the hypergraph $\lambda K_{n}^{w}$ is just $K_{n}^{w}$. The following is a generalization of Baranyai's Theorem.
\begin{Theorem}[\label{almostregular-1}\cite{Bryant}]
Let $n,w,g$ be positive integers with $\lfloor\frac{n}{w}\rfloor\mid {\lambda{{n}\choose{w}}}$.
Then hypergraph $\lambda K_{n}^{w}$ has an almost-regular edge-coloring with ${\lambda{{n}\choose{w}}}/{\lfloor\frac{n}{w}\rfloor}$ color classes such that the number of edges in each color class is $\lfloor\frac{n}{w}\rfloor$.
\end{Theorem}

Let $\mathcal{H}$ be a hypergraph.
A {\em{strong $k$-coloring}} of the vertices is a $k$-partition $S_{1},S_{2},\ldots,S_{k}$ of $V(\mathcal{H})$ such that no color appears twice in the same edge; that is to say $\mid E\cap S_{i}\mid\leq 1$ for every edge $E\in E(\mathcal{H})$, $1\leq i \leq k$. 

A {\em{cycle}} of length $k$ in a hypergraph $\mathcal{H}$ is a sequence $x_{1},E_{1},x_{2},E_{2},x_{3},\ldots,x_{k},E_{k},x_{1}$ satisfying that:\vspace{-0.2cm}
\begin{itemize}
 \item [$(1)$] $E_{1},E_{2},\ldots,E_{k}$ are distinct edges of $\mathcal{H}$;\vspace{-0.2cm}
 \item [$(2)$] $x_{1},x_{2},\ldots,x_{k}$ are distinct vertices of $\mathcal{H}$;\vspace{-0.2cm}
 \item [$(3)$] $\{x_{i},x_{i+1}\}\subseteq E_{i}$ for $1\leq i \leq k-1$ and $\{x_{k},x_{1}\}\subseteq E_{k}$.\vspace{-0.2cm}
\end{itemize}
A cycle of odd length is called an odd cycle.

By combining {\cite[pp.164, Theorem 5]{Hypergraphs}} and {\cite[pp.170, Corollary 2]{Hypergraphs}}, we obtain the following theorem.

\begin{Theorem}[\label{strong-color}\cite{Hypergraphs}]
The vertices of an $r$-uniform hypergraph without odd cycles can be strongly colored with $r$ colors.
\end{Theorem}

In the remaining part of this paper, we assume always $q\geq2$ to be the alphabet size and denote $g:=q-1$. Also let $n,w,d$  denote the length, weight and distance of a constant-weight code. Hence, $w\leq n$ and $d\le 2w$.

\section{Tilings for general weights}

In this section, we investigate several general methods for constructing $\mathrm{TOC}_{q}(n,d,w)$s, for which we distinguish them with respect to their distance $d$. When $d$ is odd, we investigate
large sets of generalized Steiner systems $\mathrm{LGS}(t,w,n,g)$, which gives a $\mathrm{TOC}_{g+1}(n,2(w-t)+1,w)$.
For the case $d\in \{2,2w\}$, we completely resolve the existence problem for $\mathrm{TOC}_{q}(n,d,w)$s.
When $d$ is even with $d\not\in \{2,2w\}$, we define large sets of generalized maximum H-packings $\mathrm{LGMHP}(t,w,n,g)$, which yields a $\mathrm{TOC}_{g+1}(n,2(w-t+1),w)$.

\subsection{Odd distances}
For $q>2$, a generalized Steiner system $\mathrm{GS}(t,w,n,g)$ is an optimal $(n,2(w-t)+1,w)_{q}$-code ($q=g+1$) by Lemmas \ref{classic-bound-generalized} and \ref{Boundgs}, so it follows that we  can consider a $\mathrm{TOC}_{q}(n,d,w)$ for  odd $d=2(w-t)+1$ such that each tile is a $\mathrm{GS}(t,w,n,g)$, for which we also name it a large set of generalized Steiner system and use the notation $\mathrm{LGS}(t,w,n,g)$.
In this subsection,
we present a recursive construction to produce $\mathrm{LGS}(t,w,n,g)$s via $t$-resolvable $\mathrm{S}(w,k,n)$s.
For $d=2w-1$, we also present  constructions to produce $\mathrm{TOC}_{q}(n,2w-1,w)$s via almost-regular edge-colorings of $K_{n}^{w}$ with certain properties.

\begin{Lemma}\label{LGS-resolvable1}
If there exists an $\mathrm{S}(t,w',n)$ and a $\mathrm{GS}(t,w,w',g)$, then there exists a $\mathrm{GS}(t,w,n,g)$.
\end{Lemma}
\proof Suppose that $([n],\mathcal{B})$ is the given $\mathrm{S}(t,w',n)$. Let $\mathcal{G}_{B}=\{\{x\}\times [g]:x\in B\}$, where $B\subseteq [n]$. For any $B\in\mathcal{B}$, we can construct by assumption a $\mathrm{GS}(t,w,w',g)$ on $B\times [g]$ with block set $\mathcal{C}_{B}$ and group set $\mathcal{G}_{B}$. Define
$$\mathcal{C} = \bigcup_{B\in \mathcal{B}}\mathcal{C}_{B}.$$
For any $t$-subset $T=\{(x_{1},a_{1}),(x_{2},a_{2}),\ldots, (x_{t},a_{t})\}\subseteq [n]\times [g]$, where  $x_{1},x_{2},\ldots, x_{t}$ are pairwise distinct, we can find exactly one block $B\in \mathcal{B}$ such that $\{x_{1},x_{2},\ldots, x_{t}\}\subseteq B$ since $([n],\mathcal{B})$ is an $\mathrm{S}(t,w',n)$. Then there is one block $C\in \mathcal{C}_{B}$ such that $T\subseteq C$ as $(B\times [g],\mathcal{G}_{B},\mathcal{C}_{B})$ is a $\mathrm{GS}(t,w,w',g)$.
So $([n]\times [g],\mathcal{G}_{[n]},\mathcal{C})$ is an $\mathrm{H}(n,g,w,t)$.

For any $C_{1}\neq C_{2}\in \mathcal{C}$, if there is a block $B\in\mathcal{B}$ such that $C_{1},C_{2}\in \mathcal{C}_{B}$, then $d_{H}(C_{1},C_{2})\geq 2(w-t)+1$ since $\mathcal{C}_{B}$ is the block set of a $\mathrm{GS}(t,w,w',g)$. If $C_{1}\in \mathcal{C}_{B_{1}}$ and $C_{2}\in \mathcal{C}_{B_{2}}$ with $B_{1}\neq B_{2}$, Then $d_{H}(C_{1},C_{2})\geq d_{H}(B_{1},B_{2})\geq 2(w'-t)+2$. It follows that $([n]\times[g],\mathcal{G}_{[n]},\mathcal{C})$ is a $\mathrm{GS}(t,w,n,g)$.
\qed

%

\begin{Theorem}\label{LGS-reslovable2}
Suppose that there exists a $t$-resolvable $\mathrm{S}(w,w',n)$. If there exists an $\mathrm{LGS}(t,w,w',g)$, then there exists an $\mathrm{LGS}(t,w,n,g)$.
\end{Theorem}
\proof Suppose $([n],\mathcal{B})$ is the given $t$-resolvable $\mathrm{S}(w,w',n)$. Thus the block set $\mathcal{B}$ can be partitioned into block sets of $\mathrm{S}(t,w',n)$, say $\mathcal{B}_{1},\mathcal{B}_{2},\ldots, \mathcal{B}_{p}$, where $p=\frac{{n\choose w}}{{w'\choose w}}\bigg/\frac{{n\choose t}}{{{w'}\choose t}}={{{n-t}\choose {w-t}}}/{{{w'-t}\choose{w-t}}}$. Let $\mathcal{G}_{B}=\{\{x\}\times [g]:x\in B\}$, where $B\subseteq [n]$.
Fix any $i$ with $1\leq i \leq p$.
For any $B\in \mathcal{B}_{i}$, we can construct an $\mathrm{LGS}(t,w,w',g)$,
say $(B\times [g],\mathcal{G}_{B},\mathcal{C}_{B,i}^{j})$, $1\leq j \leq g^{w-t}{{w'-t}\choose{w-t}}$.
For $1\leq i \leq p$ and $1\leq j \leq g^{w-t}{{w'-t}\choose{w-t}}$, define
$$\mathcal{C}_{i}^{j} = \bigcup_{B\in \mathcal{B}_{i}}\mathcal{C}_{B,i}^{j}.$$
Then each $([n]\times [g],\mathcal{G}_{[n]},\mathcal{C}_{i}^{j})$ is a $\mathrm{GS}(t,w,n,g)$ by Lemma \ref{LGS-resolvable1}.

For each $w$-subset $C=\{(x_{1},a_{1}),(x_{2},a_{2}),\ldots,(x_{w},a_{w})\}\in{\cal H}_q(n,w)$, we can find a unique block $B\in \mathcal {B}_{i}\subseteq\mathcal{B}$ such that $\{x_{1},x_{2},\ldots, x_{w}\}\subseteq B$ since $([n],\mathcal{B})$ is an $\mathrm{S}(w,w',n)$.
Next we can find some $j$ such that $C\in \mathcal{C}_{B,i}^{j}$ as the collection $\{(B\times [g],\mathcal{G}_{B},\mathcal{C}_{B,i}^{j}):1\leq j \leq g^{w-t}{{w'-t}\choose{w-t}}\}$ is an $\mathrm{LGS}(t,w,w',g)$. So the collection $$\left\{\left([n]\times [g],\mathcal{G}_{[n]},\mathcal{C}_{i}^{j}\right):1\leq i \leq p, 1 \leq j\leq g^{w-t}{{w'-t}\choose{w-t}}\right\}$$ is the desired $\mathrm{LGS}(t,w,n,g)$.
\qed

We need to introduce
a type of special almost-regular edge-coloring of the complete $w$-uniform hypergraph $K_{n}^{w}$.
 Suppose that $\lfloor\frac{n}{w}\rfloor\mid{n\choose w}$ and denote $r:={n\choose w}/\lfloor\frac{n}{w}\rfloor$. Let $g$ be positive integer with  $g\mid r$.
By Theorem \ref{almostregular}, there is an almost-regular edge-coloring  $\mathcal{F}:=\{\mathcal{F}_{i,j}:1\leq i  \leq r/g,1\leq j \leq g\}$ such that each $\mathcal{F}_{i,j}$ is an almost-regular sub-hypergraph of $K_{n}^{w}$ with $\lfloor\frac{n}{w}\rfloor$ edges. Then $\mathcal{F}$ is said to be
 {\em{$g$-good}} if  each $\mathcal{A}_{i} = \bigcup_{j=1}^{g}\mathcal{F}_{i,j}$, $1\leq i \leq r/g$, is the block set of a 2-packing $\mathrm{P}(2,w,n)$.
 Further if the vertices of each $\mathcal{A}_{i}$, $1\leq i \leq r/g$, can be strongly colored with $w$ colors, such a $g$-good almost-regular edge-coloring  is called {\em $g^{*}$-good}.

\begin{Theorem}\label{LGS-13n2-1}
Let $w,n$ be positive integers such that $r:={n\choose w}/\lfloor\frac{n}{w}\rfloor$ is even. Let $2a< w$ where $a$ is the remainder of $n$ divided by $w$. Suppose that there exists a $2$-good almost-regular edge-coloring of $K_{n}^{w}$. Then there exists a $\mathrm{TOC}_{3}(n,2w-1,w)$.
\end{Theorem}

\proof 
Suppose that there exists a $2$-good almost-regular edge-coloring of $K_{n}^{w}$ over $[n]$, say $\{\mathcal{F}_{s,t}:1\leq s \leq  r/2, t=1,2\}$, where each $\mathcal{F}_{s,t}$ forms an almost-regular sub-hypergraph of $K_{n}^{w}$ with $\lfloor\frac{n}{w}\rfloor$ edges, and each $\mathcal{A}_{s}=\mathcal{F}_{s,1}\cup\mathcal{F}_{s,2}$ with $1\leq s\leq  r/2$ is the block set of a $\mathrm{P}(2,w,n)$.
It follows that every vertex in each $\mathcal{F}_{s,t}$ has degree at most one and thus
the hypergraph $\mathcal{A}_{s}$ does not contain any odd cycle.
Thus, $\mathcal{A}_{s}$ can be strongly colored with $w$
colors by Theorem \ref{strong-color}. Correspondingly the vertex set $[n]$ is partitioned into $P_{1}^{s}, P_{2}^{s},\ldots,P_{w}^{s}$ such that $|A\cap P_{c}^{s}|=1$ for all $A\in \mathcal{A}_{s}$ and $1 \leq c  \leq w$.

Let $\mathcal{M} =(C_{1}\text{ } C_{2} \text{ }\cdots C_{w})$ be an $\mathrm{OA}(w,w,2)$ over $[2]$, where $C_{i}$ is a column vector of length $2^{w}$.
Such an $\mathrm{OA}$ exists trivially consisting of all vectors of length $w$ over $[2]$ as rows.
For $1\leq s \leq  r/2$, define a $2^{w}\times n$ array
$$\mathcal{M}_{s}=(C_{\delta_{s}(1)}\text{ }C_{\delta_{s}(2)}\text{ } \cdots C_{\delta_{s}(n)}),\ \text{where }\delta_{s}(x) =c \text{ if } x\in P_{c}^{s}.$$
For any block $\{x_{1},x_{2},\ldots,x_{w}\}\in \mathcal{A}_{s}$, letting $x_{1}<x_{2}<\cdots<x_{w}$, the sub-array $(C_{\delta_{s}(x_{1})}\ C_{\delta_{s}(x_{2})}$ $\ldots C_{\delta_{s}(x_{w})})$ of $\mathcal{M}_{s}$ is an $\mathrm{OA}(w,w,2)$ because of the strong coloring of $\mathcal{A}_{s}$.

Now, we denote by $m_{i,j}^{(s)}$ the element at the $i$-th row and $j$-th column of the array $\mathcal{M}_{s}$.
 We construct optimal $(n,2w-1,w)_{3}$-codes. For any $1\leq i \leq 2^{w}$, define
\begin{align*}
\mathcal{B}_{i}^{s}=\left\{\big\{(x_{1},m_{i,x_{1}}^{(s)}),(x_{2},m_{i,x_{2}}^{(s)}),\ldots,(x_{w},m_{i,x_{w}}^{(s)})\big\}:\{x_{1},x_{2},\ldots,x_{w}\}\in \mathcal{F}_{s,1}\right\}\\
\bigcup\left\{\big\{(y_{1},p_{i,y_{1}}^{(s)}),(y_{2},p_{i,y_{2}}^{(s)}),\ldots,(y_{w},p_{i,y_{w}}^{(s)})\big\}:\{y_{1},y_{2},\ldots,y_{w}\}\in \mathcal{F}_{s,2}\right\},
\end{align*}
where $p_{i,j}^{(s)}=1$ if $m_{i,j}^{(s)}=2$ and $p_{i,j}^{(s)}=2$ if $m_{i,j}^{(s)}=1$.
Because $\mathcal{F}_{s,1}\cup\mathcal{F}_{s,2}$ is the block set of a $\mathrm{P}(2,w,n)$, the union of the supports of any two distinct codewords in $\mathcal{B}_{i}^{s}$ has size no less than $2w-1$.
Moreover, the nonzero symbols $1$ and $2$ occur at most once in any coordinate position $x\in [n]$,
from which it follows that the distance $d\geq 2w-1$.
Because $2a< w$, each $\mathcal{B}_{i}^{s}$ has $2\lfloor\frac{n}{w}\rfloor=\lfloor\frac{2n}{w}\rfloor$ codewords.
Hence, each $\mathcal{B}_{i}^{s}$ is an optimal $(n,2w-1,w)_{3}$-code
as $A_{q}(n,2w-1,w)\leq \frac{(q-1)n}{w}$ by Lemma \ref{classic-bound-generalized}.

Now, we prove that every $w$-weight word $\{(x_{1},a_{1}),(x_{2},a_{2}),\ldots,(x_{w},a_{w})\}\in \mathcal{H}_{3}(n,w)$ occurs exactly once in the collection $\{\mathcal{B}_{i}^{s}:1\leq s\leq  r/2,1\leq i \leq 2^{w}\}$.
We can find $s$ and $t$ such that $\{x_{1},x_{2},\ldots,x_{w}\}\in \mathcal{F}_{s,t}$ as $\{\mathcal{F}_{s,t}:1\leq s \leq  r/2,1\leq t \leq 2\}$ contains all edges of  $K_{n}^{w}$ over $[n]$.
If $t=1$, we can find $i$ such that $m_{i,x_{1}}^{(s)}=a_{1}$, $m_{i,x_{2}}^{(s)}=a_{2},\ldots,m_{i,x_{w}}^{(s)}=a_{w}$ since $(M_{x_{1}}^{(s)}\ M_{x_{2}}^{(s)}\ldots M_{x_{w}}^{(s)})$ is an $\mathrm{OA}(w,w,2)$, where $M_{x_{i}}^{(s)}$ is  the $x_{i}$-th column of  $\mathcal{M}_{s}$.
Then $\{(x_{1},a_{1}),(x_{2},a_{2}),\ldots,(x_{w},a_{w})\}\in\mathcal{B}_{i}^{s}$.
If $t=2$, we also can find $i$ such that $p_{i,x_{1}}^{(s)}=a_{1}$, $p_{i,x_{2}}^{(s)}=a_{2},\ldots,p_{i,x_{w}}^{(s)}=a_{w}$,
then $\{(x_{1},a_{1}),(x_{2},a_{2}),\ldots,(x_{w},a_{w})\}\in\mathcal{B}_{i}^{s}$.
Therefore, the collection $\{\mathcal{B}_{i}^{s}:1\leq s \leq  r/2,1\leq i \leq 2^{w}\}$ is the desired $\mathrm{TOC}_{3}(n,2w-1,w)$.
\qed

The following is a generalization of Theorem \ref{LGS-13n2-1} from $q=3$ to general alphabet size.

\begin{Theorem}\label{LGS-13n2-11}
Let $w,n,g$ be positive integers such that 
$ga< w$ where $a$ is the remainder of $n$ divided by $w$. Further let both $r:={n\choose w}/\lfloor\frac{n}{w}\rfloor$ and $p:=r/g$ be integers.
Suppose that there exists a $g^{*}$-good almost-regular edge-coloring of $K_{n}^{w}$ and there exists an $\mathrm{OA}(w,2w-1,g)$. Then there exists a $\mathrm{TOC}_{q}(n,2w-1,w)$ (where $q=g+1$).
\end{Theorem}
\proof
Suppose the collection $\{\mathcal{F}_{s,t}:1\leq s \leq p, 1 \leq t\leq g\}$ is the given
$g^{*}$-good almost-regular edge-coloring of $K_{n}^{w}$ over $[n]$, where every $\mathcal{F}_{s,t}$ is an almost-regular sub-hypergraph of $K_{n}^{w}$ with $\lfloor\frac{n}{w}\rfloor$ edges.
For $1\leq s\leq p$, each $\mathcal{A}_{s}=\bigcup_{1 \leq t \leq g}\mathcal{F}_{s,t}$ is the block set of a $\mathrm{P}(2,w,n)$
and $\mathcal{A}_{s}$ can be strongly colored with $w$ colors.
Hence, the vertex set $[n]$ can be partitioned into $P_{1}^{s}, P_{2}^{s},\ldots,P_{w}^{s}$ such that $|A\cap P_{c}^{s}|=1$ for all $A\in \mathcal{A}_{s}$ and $1 \leq c  \leq w$.

Let $\mathcal{M}$ be an $\mathrm{OA}(w,2w-1,g)$ over $[g]$. Index the rows of $\mathcal{M}$ by $[g]^{w}=[g]^{w-1}\times[g]$ and then $\mathcal{M}$ is partitioned into $g$ row blocks $R_{\alpha}$, $\alpha\in [g]^{w-1}$, each of shape $g\times(2w-1)$. Since $\mathcal{M}$ is an orthogonal array, for each $\alpha=(i_{1},i_{2},\ldots,i_{w-1})\in [g]^{w-1}$, we may assume that the restriction of $R_{\alpha}$ to the last $w-1$ columns has the same $g$ rows $(i_{1}\ i_{2}\ldots i_{w-1})$ (and hence each of the first $w$ columns forms a permutations on $[g]$).

We may also denote the orthogonal array $\mathcal{M}=(C_{1}\ C_{2}\ldots C_{2w-1})$, where $C_{i}$ is a column vector of length $g^{w}$. For $1\leq s \leq p$, define a $g^{w}\times n$ array
$$\mathcal{M}_{s}=(C_{\delta_{s}(1)}\ C_{\delta_{s}(2)}\ldots C_{\delta_{s}(n)}),\  \text{where }\delta_{s}(x)=c \text{ if }x\in P_{c}^{s}.$$
For any block $\{x_{1},x_{2},\ldots,x_{w}\}\in \mathcal{A}_{s}$, letting $x_{1}<x_{2}<\cdots<x_{w}$, the sub-array $(C_{\delta_{s}(x_{1})}\ C_{\delta_{s}(x_{2})}$ $\ldots  C_{\delta_{s}(x_{w})})$ is an $\mathrm{OA}(w,w,g)$, due to the strong coloring of $\mathcal{A}_{s}$. We also index the rows of $\mathcal{M}_{s}$ by $[g]^{w-1}\times[g]$ and denote $m_{(\alpha,i),j}^{(s)}$ its entry in the row $(\alpha,i)$ and the column $j$, where $\alpha\in [g]^{w-1}$, $i\in [g]$, $j\in [n]$.

Now, we construct optimal $(n,2w-1,w)_{q}$-codes. Note that $A_{q}(n,2w-1,w)\leq \lfloor\frac{gn}{w}\rfloor=g\lfloor\frac{n}{w}\rfloor$ if $ag<w$. For $ s \in [p]$, $t,i\in [g]$, $\alpha\in [g]^{w-1}$, define
$$\mathcal{C}_{s,t}(\alpha,i)=\bigg\{\left\{(x_{1},m_{(\alpha,i),x_{1}}^{(s)}),(x_{2},m_{(\alpha,i),x_{2}}^{(s)}),\ldots, (x_{w},m_{(\alpha,i),x_{w}}^{(s)})\right\}:\{x_{1},x_{2},\ldots,x_{w}\}\in\mathcal{F}_{s,t}\bigg\}$$
and
$$\mathcal{C}_{s}(\alpha,i)=\bigcup_{1\leq t \leq g}\mathcal{C}_{s,t}(\alpha,i+t),$$
where $i+t$ is reduced modulo $g$ to lie in $[g]$.

Fix $s$ with $s \in [p]$ and $\alpha,i$ with $\alpha\in [g]^{w-1}$, $i\in [g]$. Note that
$\mathcal{F}_{s,1},\mathcal{F}_{s,2},\ldots,\mathcal{F}_{s,g}$ are disjoint almost-regular sub-hypergraphs of $K_{n}^{w}$ with $\lfloor\frac{n}{w}\rfloor$ edges and the union of $\mathcal{F}_{s,1},\mathcal{F}_{s,2},\ldots,$ $\mathcal{F}_{s,g}$ is the block set of a $\mathrm{P}(2,w,n)$. Then it follows that the intersection of the supports of any two codewords of $\mathcal{C}_{s}(\alpha,i)$ has size at most one.
Because the set $\{i+t:1\leq t \leq g\}=[g]$, where $i+t\in [g]$ modulo $g$,
each nonzero symbol of $[g]$ occurs at most once in any given coordinate position.
Then the  distance of $\mathcal{C}_{s}(\alpha,i)$ is $2w-1$.
Because
$\mathcal{C}_{s}(\alpha,i)$ has $\lfloor\frac{gn}{w}\rfloor$ codewords,
 each $\mathcal{C}_{s}(\alpha,i)$ is an optimal $(n,2w-1,w)_{q}$-code.

Finally, we prove that every $w$-weight word $\{(x_{1},a_{1}),(x_{2},a_{2}),\ldots,(x_{w},a_{w})\}\in \mathcal{H}_{q}(n,w)$ occurs exactly once in the collection $\{\mathcal{C}_{s}(\alpha,i):s \in [p],\alpha\in [g]^{w-1},i\in [g]\}$.
We can find $s$ and $t$ such that $\{x_{1},x_{2},\ldots,x_{w}\}\in \mathcal{F}_{s,t}$ as $\{\mathcal{F}_{s,t}:1\leq s \leq p,1\leq t \leq g\}$ is a partition of edges in $K_{n}^{w}$ over $[n]$.
We can find $\alpha\in [g]^{w-1},i\in [g]$ such that $\{(x_{1},a_{1}),(x_{2},a_{2}),\ldots,(x_{w},a_{w})\}\in\mathcal{C}_{s,t}(\alpha,i)$ since the sub-array $(C_{\delta_{s}(x_{1})},C_{\delta_{s}(x_{2})},$ $\ldots, C_{\delta_{s}(x_{w})})$ is an $\mathrm{OA}(w,w,g)$.
So, we obtain a $\mathrm{TOC}_{q}(n,2w-1,w)$ indeed by noting Lemma \ref{mq1}.
\qed

%

\subsection{Even distances}
This subsection starts with the cases $d=2$ and $d=2w$, for which a $\mathrm{TOC}_{q}(n,d,w)$ is not difficult to be constructed.
When $d$ is even with $d=2(w-t+1)$ and $d\notin\{2,2w\}$, we define large sets of generalized maximum $\mathrm{H}$-packings, as an analog of large sets of generalized
Steiner systems.
Then we develop some analogous constructions.

\begin{Lemma}[\label{bound_A_q_n2w}\cite{bound_d2}]
$A_{q}(n,2,w)={{n}\choose{w}}(q-1)^{w-1}$, $A_{q}(n,2w,w)=\lfloor\frac{n}{w}\rfloor.$
\end{Lemma}

\begin{Theorem}\label{TOC_n_2_w}
There exists a $\mathrm{TOC}_{q}(n,2,w)$ for any $q,n$ and $w\leq n$.
\end{Theorem}
\proof
For any $1\leq x_{1}<\cdots<x_{w}\leq n$ and $1\leq i \leq g=q-1$, let
$$\mathcal{C}_{i}(x_{1},\ldots,x_{w})=\big\{\left\{(x_{1},a),(x_{2},a+i),(x_{3},b_{3}),\ldots,(x_{w},b_{w})\right\}:a,b_{3},\ldots,b_{w}\in [g]\big\},$$
where $a+i$ is reduced modulo $g$ to lie in $[g]$.
In each $\mathcal{C}_{i}(x_{1},\ldots,x_{w})$, the support of each codeword is $\{x_{1},x_{2},\ldots,x_{w}\}$ and
the nonzero symbols in the coordinate positions $x_{1}$ and $x_{2}$ are different for any two different codewords in $\mathcal{C}_{i}(x_{1},\ldots,x_{w})$.
Then the minimum distance of $\mathcal{C}_{i}(x_{1},\ldots,x_{w})$ is $2$.
As a result, for $1 \leq i\leq g$, the code
$$\mathcal{C}_{i}=\bigcup_{1\leq x_{1}<\ldots<x_{w}\leq n}\mathcal{C}_{i}(x_{1},\cdots,x_{w})$$
is obviously an optimal $(n,2,w)_{q}$-code with $g^{w-1}{{n}\choose{w}}$ codewords.

For any $w$-weight word $\{(x_{1},b_{1}),(x_{2},b_{2}),\ldots, (x_{w},b_{w})\}\in \mathcal{H}_{q}(n,w)$ with $x_{1}<x_{2}<\cdots<x_{w}$, let $b_{2}-b_{1}\equiv i\pmod{g}$ with $i\in [g]$. Then we have $\{(x_{1},b_{1}),(x_{2},b_{2}),\ldots, (x_{w},b_{w})\}\in \mathcal{C}_{i}(x_{1},\ldots,x_{w}) \subseteq \mathcal{C}_{i}$.
Because $g\cdot g^{w-1}{n\choose w}=|\mathcal{H}_{q}(n,w)|$, the collection $\{\mathcal{C}_{i}:1\leq i \leq g\}$
 is a $\mathrm{TOC}_{q}(n,2,w)$.
\qed

When $d=2w$,  the supports of codewords in an  $(n,2w,w)$-code are mutually disjoint. 
Applying Theorem \ref{almostregular-1}, we have the following theorem.


\begin{Theorem}\label{TOC_n_2w_w_1}
For any $q,n$ and $w\leq n$, there exists a $\mathrm{TOC}_{q}(n,2w,w)$ if and only if $\lfloor\frac{n}{w}\rfloor$ divides ${(q-1)^{w}{{n}\choose{w}}}$.
\end{Theorem}
\proof The necessity is obvious by noting Lemma \ref{bound_A_q_n2w}. For sufficiency, let $\lambda=g^{w}$ (again $g=q-1$). By Theorem \ref{almostregular-1},
if $\lfloor\frac{n}{w}\rfloor\mid{\lambda{{n}\choose{w}}}$,
the hypergraph $\lambda K_{n}^{w}$ has an almost-regular edge-coloring  with color classes  $\mathcal{F}_{1},\mathcal{F}_{2},\ldots,\mathcal{F}_{k}$,
where $k={g^{w}{{n}\choose{w}}}/{\lfloor\frac{n}{w}\rfloor}$ and each $\mathcal{F}_{i}$, $1\leq i \leq k$, has $\lfloor\frac{n}{w}\rfloor$ pairwise disjoint edges.
Note that 
each $w$-subset $A$ of $[n]$  appears $\lambda$ times in the hypergraph $\lambda K_{n}^{w}$. Thus $A$ is contained in $\lambda=g^{w}$ color classes, say $\mathcal{F}_{i_1},\mathcal{F}_{i_2},\ldots,\mathcal{F}_{i_\lambda}$. Arbitrarily replace $A$ in every  ${\cal F}_{i_j}$ with one of the $g^w$ words  in $\mathcal{H}_{q}(n,w)$ with support $A$, but we need run out of all such words eventually. Repeat this process for every  $w$-subset. Then we produce a code $\mathcal{C}_{i}$ of weight $w$ from the  color class $\mathcal{F}_{i}$ for every  $1\leq i \leq k$. Each code $\mathcal{C}_{i}$ has  $\lfloor\frac{n}{w}\rfloor$ codewords with distance $2w$ because the supports  are pairwise disjoint. It is immediate that the collection $\big\{\mathcal{C}_{i}:1\leq i \leq k\big\}$ is a $\mathrm{TOC}_{q}(n,2w,w)$.
\qed

%

When $d$ is even with $d\notin \{2,2w\}$, let $t=\frac{2w-d+2}{2}$.
In Subsection 2.2, we showed that a  $\mathrm{GMHP}(t,w,n,g)$ is an optimal  $(n,d,w)_q$-code with  even distance $d=2(w-t+1)$.  In this subsection we  consider $\mathrm{TOC}_{q}(n,d,w)$ with $d=2(w-t+1)$ such that each tile is a $\mathrm{GMHP}(t,w,n,g)$. For such a tiling, we also use the name large set and denote it $\mathrm{LGMHP}(t,w,n,g)$. By Lemma \ref{classic-bound-generalized}, we have $A_{q}(n,d,w)\leq g^{t-1}{{n}\choose{t}}/{{w}\choose{t}}$. Observing that this upper bound is not tight always, we attach an asterisk to $\mathrm{GMHP}(t,w,n,g)$
if it contains $ g^{t-1}{{n}\choose{t}}/{{w}\choose{t}}$ blocks.
By applying construction methods analogous to that used in Lemma \ref{LGS-resolvable1} and Theorem \ref{LGS-reslovable2}, we establish the following.

\begin{Lemma}\label{LGS-resolvable3}
If there exists an $\mathrm{S}(t,w',n)$ and a $\mathrm{GMHP}^*(t,w,w',g)$, then there exists a $\mathrm{GMHP}^*$ $(t,w,n,g)$.
\end{Lemma}

\proof
It is similar to the proof of Lemma \ref{LGS-resolvable1}. Start from an  $\mathrm{S}(t,w',n)$ and for each of its blocks construct a  $\mathrm{GMHP}^*(t,w,w',g)$ instead of  $\mathrm{GS}(t,w,w',g)$ to obtain the desired  $\mathrm{GMHP}^*$$(t,w,n,g)$.
We omit the detailed description here and only check the count. An $\mathrm{S}(t,w',n)$ has $\frac{{n\choose t}}{{w'\choose t}}$ blocks and a $\mathrm{GMHP}^*(t,w,w',g)$ has $\frac{{w'\choose t}g^{t-1}}{{w\choose t}}$ blocks. Then we get a $\mathrm{GMHP}(t,w,n,g)$ with
$$\frac{{n\choose t}}{{w'\choose t}}\cdot \frac{{w'\choose t}g^{t-1}}{{w\choose t}}=\frac{{n\choose t}g^{t-1}}{{w\choose t}} $$ blocks, that is,
 a $\mathrm{GMHP}^*(t,w,n,g)$.\qed

\begin{Theorem}\label{LGS-reslovable4}
Suppose that there exists a $t$-resolvable $\mathrm{S}(w,w',n)$. If there exists  an $\mathrm{LGMHP}^*$ $(t,w,w',g)$, then there exists an $\mathrm{LGMHP}^*(t,w,n,g)$.
\end{Theorem}
\proof It is a complete analogy to the proof of Theorem \ref{LGS-reslovable2}, also noting $$m_q(n,2(w-t+1),w)=g^{w-t+1}{{n-t}\choose{w-t}}
={{n-t}\choose {w-t}}\bigg/{{w'-t}\choose{w-t}}\times  m_q(w',2(w-t+1),w).$$
\qed
%

\section{Tilings for weight three}

This section  pays attention to $\mathrm{TOC}_{q}(n,d,3)$. For binary case,
the existence problem of  $\mathrm{TOC}_{2}(n,d,3)$ can be totally resolved.  For small alphabet size $q\ge 3$, we obtain many infinite families of $\mathrm{TOC}_{q}(n,d,3)$ for distances $d=3,4,5$.


\subsection{$q=2$}

It is immediate that $A_{2}(n,d-1,w)=A_{2}(n,d,w)$ for any even $d$ with $d\leq 2w$.
The following results are derived from optimal packings $\mathrm{P}(2,3,n)$ {\cite{bound_q2}}.
\begin{Lemma}[\label{Boundq2}\cite{bound_q2}]
\begin{itemize}
\item [$(1)$] If $n\not\equiv 5\pmod{6}$, then $A_{2}(n,4,3)=\lfloor\frac{n}{3}\lfloor\frac{n-1}{2}\rfloor\rfloor$.\vspace{-0.2cm}
\item [$(2)$] If $n\equiv 5\pmod{6}$, then
$A_{2}(n,4,3)=\lfloor\frac{n}{3}\lfloor\frac{n-1}{2}\rfloor\rfloor-1$.
\end{itemize}
\end{Lemma}


\begin{Theorem}\label{Boundq2--T} Consider tilings with $q=2$ and $w=3$. We have the following:
\begin{itemize}
\item [$(1)$] For $d=1,2$, a $\mathrm{TOC}_{2}(n,d,3)$ exists for all $n\geq3$;\vspace{-0.2cm}
\item [$(2)$] For $d=3,4$, a $\mathrm{TOC}_{2}(n,d,3)$ exists if and only if $n=3,4,5$ or $n \equiv 0,1,2,3\pmod{6}$, $n\geq  8$;\vspace{-0.2cm}
\item [$(3)$] For $d=5,6$, a $\mathrm{TOC}_{2}(n,d,3)$ exists if and only if $n \not\equiv 1\pmod{6}$ and $n\geq3$.
\end{itemize}
\end{Theorem}
\proof Part $(1)$ follows from Theorem  \ref{TOC_n_2_w}. If $n\geq3$ and $n \not\equiv 1\pmod{6}$, a $\mathrm{TOC}_{2}(n,6,3)$ exists by Theorem  \ref{TOC_n_2w_w_1}. If $n \equiv 1\pmod{6}$,  a $\mathrm{TOC}_{2}(n,6,3)$ does not exist because $A_2(n,6,3)={n-1\over 3}$ does not divide $|{\cal H}_2(n,3)|= {n\choose 3}$. This proves Part $(3)$.

Now we prove (2). For $n=3,4$, the existence of  a $\mathrm{TOC}_{2}(n,4,3)$ is trivial, because the corresponding optimal CWC contains only one codeword.
For $n=5$,  a $\mathrm{TOC}_{2}(5,4,3)$ exists with five codes: $\{(11100),(00111)\},$ $\{(11010),(01101)\}$, $ \{(01110), (10011)\},$ $\{(11001),(10110)\},$ $\{(10101),(01011)\}.$
For $n\equiv 4,5\pmod{6}$ and $n\geq 10$,  it can be checked that $A_{2}(n,4,3)\nmid {n\choose 3}$ by Lemma \ref{Boundq2}.
Therefore, we only need to consider $n \equiv 0,1,2,3\pmod{6}$. Since an optimal $(n,4,3)_2$-code is equivalent to an optimal packing $\mathrm{P}(2,3,n)$, it follows that a $\mathrm{TOC}_{2}(n,4,3)$ is equivalent to a partition of all triples in $[n]$ into  optimal packings $\mathrm{P}(2,3,n)$.  Note that an optimal $\mathrm{P}(2,3,n)$
is an $\mathrm{STS}(n)$ if $n \equiv 1,3\pmod{6}$ and is a $\mathrm{P}(2,3,n)$ with $n(n-2)/6$ blocks if $n \equiv 0,2\pmod{6}$. The optimal $\mathrm{P}(2,3,n)$ for $n \equiv 0,2\pmod{6}$ can be obtained from an $\mathrm{STS}(n+1)$ by deleting a point and the vice versa (by adjoining a new point to all pairs not be contained in any block). This correspondence could be extended to the equivalence between a $\mathrm{TOC}_{2}(n,4,3)$ for  $n \equiv 1,3\pmod{6}$ and a $\mathrm{TOC}_{2}(n,4,3)$ for  $n \equiv 0,2\pmod{6}$.
It is well-known that a large set of $\mathrm{STS}(n)$s exists if and only if
$n \equiv 1,3\pmod{6}$, $n\geq 3$ and $n\neq 7$, which directs us to the conclusion in Part (2).
%
%
%
\qed

\subsection{$d=3$}

As an application of  Theorem \ref{LGS-reslovable2}, we display the existence of infinite families of $\mathrm{TOC}_{3}(n,3,3)$s.

\begin{Example}\label{LGS2342}
{\rm The following lists a $\mathrm{TOC}_{3}(4,3,3)$.}\vspace{-0.2cm} 
\end{Example}\vspace{-0.2cm}
\begin{center}\vspace{-0.2cm}
\begin{tabular}{l l l l l l l l }
$\{(1\text{ } 1\text{ }1 \text{ }0),(2\text{ } 2\text{ } 2\text{ }0),(2\text{ } 1\text{ } 0\text{ }2),(1\text{ } 2\text{ } 0\text{ }1),(1\text{ } 0\text{ } 2\text{ }2),(2\text{ } 0\text{ } 1\text{ }1),(0\text{ } 2\text{ } 1\text{ }2),(0\text{ } 1\text{ } 2\text{ }1)\}$,\\
$\{( 2\text{ }1  \text{ } 1 \text{ }0 ),( 1\text{ } 2 \text{ } 2 \text{ }0 ),( 1\text{ } 1 \text{ }  0\text{ }1 ),( 2\text{ } 2 \text{ }  0\text{ }2 ),(1 \text{ } 0 \text{ }1  \text{ } 2 ),(2 \text{ } 0 \text{ }  2\text{ } 1),( 0\text{ } 1 \text{ } 2 \text{ } 2),( 0\text{ } 2 \text{ }  1\text{ }1 )\}$,\\
$\{( 1\text{ }1  \text{ } 2 \text{ } 0),( 2\text{ } 2 \text{ } 1 \text{ } 0),( 1\text{ } 2 \text{ }0  \text{ } 2),( 2\text{ }  1\text{ } 0 \text{ } 1),( 2\text{ }0  \text{ } 2 \text{ } 2),( 1\text{ } 0 \text{ } 1 \text{ } 1),(0 \text{ } 1 \text{ } 1 \text{ }2 ),( 0\text{ } 2 \text{ }  2\text{ } 1)\}$,\\
$\{( 1\text{ } 2 \text{ }1  \text{ } 0),( 2\text{ }  1\text{ }  2\text{ } 0),( 1\text{ } 1 \text{ } 0 \text{ } 2),( 2\text{ }  2\text{ }  0\text{ }1 ),(2 \text{ }  0\text{ } 1 \text{ } 2),( 1\text{ } 0 \text{ }  2\text{ } 1),( 0\text{ } 1 \text{ }1  \text{ } 1),(0 \text{ } 2 \text{ } 2 \text{ }2 )\}$.\\
\end{tabular}
\end{center}

\begin{Theorem}\label{LGS23ng1}
Let $K = \{4^{m}:m\geq 1\}\cup\{2\cdot q^{m}+2:q =7,31,127,m\geq 1\}$.
There exists a $\mathrm{TOC}_{3}(n,3,3)$ for any $n \in K$.
\end{Theorem}
\proof There exists a $2$-resolvable $\mathrm{S}(3, 4, n)$ for $n \in K$ by Theorem \ref{2RSQS1}. Then apply Theorem \ref{LGS-reslovable2}.
 Since a $\mathrm{TOC}_{3}(4,3,3)$, namely an $\mathrm{LGS}(2,3, 4,2)$, exists by Example \ref{LGS2342},  there exists an $\mathrm{LGS}(2,3, n,2)$ or a $\mathrm{TOC}_{3}(n,3,3)$ for $n \in K$.
 \qed

\subsection{$d=4$}

This subsection makes use of orthogonal arrays to construct a $\mathrm{TOC}_{q}(n,4,3)$ of length $n=4$ for any $q$, or of larger length $n$ for $q-1$ a prime power.
 Large sets of Steiner triple systems and their generalizations are  also employed.

\begin{Lemma}[\label{A3n43}{\cite{Chee-2008}}]
 $A_{q}(n,4,3)=\min\big\{{{n}\choose{3}},U_{q}(n)\big\},$
 where
$$
U_{q}(n)=
\begin{cases}
\big\lfloor {(q-1)n\over 3}\big\lfloor {n-1\over 2}\big\rfloor\big\rfloor -1,& \text{ if } n\equiv 5\pmod 6\text{ and }q\not\equiv 1\pmod{3};\\
\big\lfloor {(q-1)n\over 3}\big\lfloor {n-1\over 2}\big\rfloor\big\rfloor,& \text{ otherwise}.\\
\end{cases}
$$
\end{Lemma}

\begin{Theorem}\label{n43_n=4}
There exists a $\mathrm{TOC}_{q}(4,4,3)$ for any alphabet size $q$.
\end{Theorem}
\proof For $q=2,3$, see Theorem \ref{Boundq2--T} and Example \ref{(4,4,3)_3}. Next let $q\ge 4$. Then $A_{q}(4,4,3)=4$ by Lemma \ref{A3n43}. The desired tiling consists of the following $g^3$ codes ${\cal C}_{x,y,z}$, $x,y,z\in [g]$:
\begin{align*}
{\cal C}_{x,y,z}=\left\{
  \begin{array}{l}
  \{(1,x),(2,y),(3,z)\},\{(1,x+1),(2,y+1),(4,w)\},
 \{(1,x+2),(3,z+1),(4,w+2)\},\\
\ \ \{(2,y+2),(3,z+2),(4,w+1)\}:
 \ x+y+z+w\equiv 0\pmod g,\  w\in [g] \end{array}
 \right\},
\end{align*}
where the addition is computed under modulo $g$ such that the result falls in  $[g]$. It is easy to know that these codes  have distance 4. To prove that any word in ${\cal H}_q(4,3)$ appears once in the tiling, we show for instance the word $\{(1,a),(3,b),(4,c)\}$. Letting $x=a-2,z=b-1$, there must exist some $y\in [g]$ such that  $x+y+z+c\equiv 0\pmod g$. Hence we could find it in the code ${\cal C}_{x,y,z}$.
\qed


\begin{Proposition}\label{n43_gen}
Let $r=\lceil{n\choose 3}/{\lfloor\frac{n}{3}\rfloor}\rceil$.
If there exists an $\mathrm{OA}(3,n+2,g)$ for $g\ge r$, then there exists a $\mathrm{TOC}_{q}(n,4,3)$ $(q=g+1)$.
\end{Proposition}
\proof
Let $r=\lceil{n\choose 3}/{\lfloor\frac{n}{3}\rfloor}\rceil$.
Suppose that $\mathcal{M}$ is the given $\mathrm{OA}(3,n+2,g)$ over $[g]$ with $g\ge r$
and $\mathcal{M}$ is partitioned into $g^2$ row blocks $M_{i,j},i,j\in [g]$, where
\begin{align}\label{OAM}
M_{i,j}=\left(\begin{matrix}
m_{i,j,1}^{(1)} & m_{i,j,1}^{(2)} &\ldots& m_{i,j,1}^{(n)} & j & i\\
m_{i,j,2}^{(1)} & m_{i,j,2}^{(2)} &\ldots& m_{i,j,2}^{(n)} & j & i\\
\vdots&\vdots&&\vdots&\vdots&\vdots\\
m_{i,j,g}^{(1)} & m_{i,j,g}^{(2)} &\ldots& m_{i,j,g}^{(n)} & j & i\\
\end{matrix}\right).
\end{align}

Because $g\ge r$, then $A_{q}(n,4,3)={n\choose 3}$ by Lemma \ref{A3n43}.
There exists an almost-regular edge-coloring of $K_{n}^{3}$ over $[n]$ with $r$ color classes by Baranyai's Theorem \cite{Baranyai}, say $\mathcal{F}_{s}$, $1\leq s \leq r$, where each $\mathcal{F}_{s}$ with $1\leq s\leq r-1$ has $\lfloor\frac{n}{3}\rfloor$ edges and $\mathcal{F}_{r}$ has  ${n\choose 3}-(r-1)\lfloor\frac{n}{3}\rfloor$ edges.
For $i,j,k,s\in [g]$, define
$$\mathcal{B}_{i,j,k}^{s}=\big\{\{(x,m_{i,j,k}^{(x)}),(y,m_{i,j,k}^{(y)}),(z,m_{i,j,k}^{(z)})\}:\{x,y,z\}\in \mathcal{F}_{s}\big\},$$
and
$$\mathcal{C}_{i,j,k}=\bigcup_{1\leq s \leq r}\mathcal{B}_{i,j,k+s}^{s},$$
where $k+s$ is reduced modulo $g$ to lie in $[g]$.
Fix $i,j,k$ with $i,j,k\in [g]$. For any two codewords $C_{1}=\{(x_{1},a_{1}),(y_{1},b_{1}),(z_{1},c_{1})\}$, $C_{2}=\{(x_{2},a_{2}),(y_{2},b_{2}),(z_{2},c_{2})\}$ in $\mathcal{C}_{i,j,k}$, we have  $|\{x_{1},y_{1},z_{1}\}\cup\{x_{2},y_{2},z_{2}\}|\geq 4$ and $C_{1}\cap C_{2}=\emptyset$. Then $d_{H}(C_{1},C_{2})\geq 4$. Now each $\mathcal{C}_{i,j,k}$ has $n\choose 3$ codewords, so it is an optimal $(n,4,3)_{q}$-code.

For any $3$-weight word $\{(x,\alpha),(y,\beta),(z,\gamma)\}\in \mathcal{H}_{q}(n,3)$, we can find $i,j,t$ such that $m_{i,j,t}^{(x)}=\alpha$, $m_{i,j,t}^{(y)}=\beta$ and $m_{i,j,t}^{(z)}=\gamma$ as $\mathcal{M}$ is an $\mathrm{OA}(3,n+2,g)$. We can find $s$ such that $\{x,y,z\}\in \mathcal{F}_{s}$ since $\mathcal{F}_{i}$, $1 \leq i \leq r$ is an almost-regular edge-coloring of $K_{n}^{3}$.
Let $k+s\equiv t\pmod{g}$ and $k\in[g]$. Then $\{(x,\alpha),(y,\beta),(z,\gamma)\}\in \mathcal{B}_{i,j,t}^{s}\subseteq\mathcal{C}_{i,j,k}$.
So the collection $\{\mathcal{C}_{i,j,k}:i,j,k\in [g]\}$ is a $\mathrm{TOC}_{q}(n,4,3)$.
\qed

Notice that in the above proposition the alphabet size is rather large. To improve this,
  we generalize large sets of Steiner triple systems, on the basis of which a new construction will be developed.
A $\mathrm{PT}(n)$ is a partition of all triples in an $n$-set into pairwise disjoint block sets $\mathcal{B}_{i}$, $1\leq i \leq m$, such that each $\mathcal{B}_{i}$ is the block set of a $\mathrm{P}(2,3,n)$.
If the number $m$ of pairwise disjoint $\mathrm{P}(2,3,n)$s achieves its minimum, we refer to $\mathrm{PT}(n)$ as optimal and denote it by $\mathrm{OPT}(n)$.
The size $m$ of an  $\mathrm{OPT}(n)$ has been determined by combining a lot of research papers.
For $n\equiv 1,3\pmod{6}$ and $n>7$ or $n=3$, a large set of $\mathrm{STS}(n)$s exists, and the number of $\mathrm{STS}(n)$s is $n-2$.
For $n\equiv 0,2\pmod{6}$ and $n>6$, a large set of $\mathrm{STS}(n+1)$s exists, say $\{([n+1],\mathcal{A}_{i}):1\leq i \leq n-1\}$. Then the block set $\mathcal{B}_{i}$ derived from $\mathcal{A}_{i}$ defined as $\{B:n+1\not\in B,B\in \mathcal{A}_{i}\}$, is an optimal $\mathrm{P}(2,3,n)$. So the collection $\{([n],\mathcal{B}_{i}):1\leq i \leq n-1\}$ is an $\mathrm{OPT}(n)$.
For $n\equiv 4\pmod{6}$, there exists an $\mathrm{OPT}(n)$ with $n-1$ optimal packings and one packing of size $\frac{n-1}{3}$ in \cite{Etzion1992}.
For $n\equiv 5\pmod{6}$, there exists an $\mathrm{OPT}(n)$ with $n-2$ optimal packings and one packing of size $\frac{4n-8}{3}$ in \cite{Ji-2006-GC}.
For $n=7$,  an $\mathrm{OPT}(7)$ exists with six block sets of $\mathrm{P}(2,3,7)$:
\vspace{-0.2cm}
\begin{center}
\begin{tabular}{l l l l l l l l }
$\{1\ 3\ 4, \ 1\ 5\ 6,\ 2\ 3\ 5,\ 2\ 4\ 6,\ 1\ 2\ 7,\ 3\ 6\ 7,\ 4\ 5\ 7\},\ $\\
$\{1\ 2\ 6,\ 1\ 3\ 5,\ 2\ 4\ 5,\ 3\ 4\ 6,\ 1\ 4\ 7,\ 2\ 3\ 7,\ 5\ 6\ 7\},\ $\\
$\{1\ 2\ 5,\ 1\ 4\ 6,\ 2\ 3\ 4,\ 3\ 5\ 6,\ 1\ 3\ 7,\ 2\ 6\ 7\},\ $\\
$\{1\ 2\ 4,\ 1\ 3\ 6,\ 2\ 5\ 6,\ 3\ 4\ 5,\ 1\ 5\ 7,\ 4\ 6\ 7\},\ $\\
$\{1\ 2\ 3,\ 1\ 4\ 5,\ 1\ 6\ 7,\ 2\ 4\ 7,\ 3\ 5\ 7\},\ $\\
$\{2\ 3\ 6,\ 4\ 5\ 6,\ 2\ 5\ 7,\ 3\ 4\ 7\}$.\\
\end{tabular}
\end{center}\vspace{-0.2cm}
For $n=6$, an $\mathrm{OPT}(6)$ can be obtained from an $\mathrm{OPT}(7)$ by deleting a fixed point.
To summarize, the size $m$ of an  $\mathrm{OPT}(n)$ has been determined as follows:\begin{equation}\label{m-opt}
m=
\begin{cases}
n-2,& \text{ if } n\equiv 1,3\pmod 6,n\ne 7;\\
n-1,& \text{ if }  n=7 \text{ or } n\equiv 0,2,5\pmod 6,n\ne 6;\\
n,& \text{ if } n=6\text{ or } n\equiv 4\pmod 6.\\
\end{cases}
\end{equation}

\begin{Theorem}\label{n43_q2}
For $n\ge 3$, 
if there exists an $\mathrm{OA}(3,n+2,g)$,
then there exists a $\mathrm{TOC}_{g+1}(n,4,3)$ if $g\geq m$, where $m$ is defined in $(\ref{m-opt})$.
\end{Theorem}
\proof As stated previously, there exists  an $\mathrm{OPT}(n)$ on the point set $[n]$ with $m$ block sets of $\mathrm{P}(2,3,n)$, say
 $\mathcal{A}_{i},$ $1\leq i \leq m$. Let $g\geq m$. It is readily checked that $A_{g+1}(n,4,3)={n\choose 3}$ by Lemma \ref{A3n43}.
Suppose that $\mathcal{M}$ is the given $\mathrm{OA}(3,n+2,g)$ over $[g]$
and $\mathcal{M}$ has $g^2$ row blocks $M_{i,j},i,j\in [g]$, of the form (\ref{OAM}).

For $i,j,k,s\in [g]$, define
$$\mathcal{B}_{i,j,k}^{s}=\big\{\{(x,m_{i,j,k}^{(x)}),(y,m_{i,j,k}^{(y)}),(z,m_{i,j,k}^{(z)})\}:\{x,y,z\}\in \mathcal{A}_{s}\big\},$$
and
$$\mathcal{C}_{i,j,k}=\bigcup_{1\leq s \leq m}\mathcal{B}_{i,j,k+s}^{s},$$
where $k+s\in [g]$ by modulo $g$.
Fix $i,j,k$ with $i,j,k\in [g]$. For any two codewords $C_{1}=\{(x_{1},a_{1}),(y_{1},b_{1}),(z_{1},c_{1})\}$, $C_{2}=\{(x_{2},a_{2}),(y_{2},b_{2}),(z_{2},c_{2})\}$ in $\mathcal{C}_{i,j,k}$, we can find $s_{1},s_{2}$ such that $\{x_{1},y_{1},z_{1}\}\in \mathcal{A}_{s_{1}}$ and $\{x_{2},y_{2},z_{2}\}\in \mathcal{A}_{s_{2}}$ as the collection $\{\mathcal{A}_{i}:1\leq i \leq m\}$ is an $\mathrm{OPT}(n)$ on the point set $[n]$.
If $s_{1}=s_{2}$, we have $|\{x_{1},y_{1},z_{1}\}\cup\{x_{2},y_{2},z_{2}\}|-|\{x_{1},y_{1},z_{1}\}\cap\{x_{2},y_{2},z_{2}\}|\geq 4$ and then $d_{H}(C_{1},C_{2})\geq 4$. If $s_{1}\neq s_{2}$,
we have  $|\{x_{1},y_{1},z_{1}\}\cup\{x_{2},y_{2},z_{2}\}|\geq 4$ and $C_{1}\cap C_{2}=\emptyset$. So $d_{H}(C_{1},C_{2})\geq 4$ holds always. Now that each $\mathcal{C}_{i,j,k}$ has $n\choose 3$ codewords, it is an optimal $(n,4,3)_{g+1}$-code.
Analogous to the proof of Proposition \ref{n43_gen}, it can be proved that the collection $\{\mathcal{C}_{i,j,k}:i,j,k\in[g]\}$ is a partition of $\mathcal{H}_{g+1}(n,3)$, and it is a $\mathrm{TOC}_{g+1}(n,4,3)$.
\qed

Let $g$ be a prime power and $g\geq n$ if $g$ is even and $g\geq n+1$ otherwise.
There exists an $\mathrm{OA}(3,n+2,g)$ by Theorem \ref{OA_B4}.
So applying Theorem \ref{n43_q2} we have the following result.

\begin{Theorem}\label{d4gprime}
Let $g$ be a prime power and $g\geq n$ if $g$ is even and $g\geq n+1$ otherwise.
There exists a $\mathrm{TOC}_{g+1}(n,4,3)$.
\end{Theorem}

\subsection{$d=5$}

This subsection  applies  Theorems \ref{LGS-reslovable2} and \ref{LGS-13n2-1} to construct infinite families of $\mathrm{TOC}_{q}(n,$ $5,3)$s for small alphabets.

In Subsection \ref{H+LS} we introduced the concepts of large sets and overlarge sets of Kirkman triple systems. Obviously, an  $\mathrm{LKTS}(n)$  gives a 1-factorization of $K_{n}^{3}$ and  an  $\mathrm{OLKTS}(n+1)$  gives a near 1-factorization of $K_{n+1}^{3}$. Furthermore, these large sets may have the 2-good property.
A $\mathrm{KTS}(n)$ has $\frac{n-1}{2}$ parallel classes.   Note that $2\mid \frac{n-1}{2}$ if $n\equiv 9 \pmod{12}$. It follows that if there exists an $\mathrm{LKTS}(n)$ for $n\equiv 9 \pmod{12}$, then there exists a $2$-good almost-regular edge-coloring of $K_{n}^{3}$.
By the same token, if there exists an $\mathrm{OLKTS}(n)$ for $n\equiv 9 \pmod{12}$, then there exists a $2$-good almost-regular edge-coloring of $K_{n+1}^{3}$.
So we have the following corollary to Theorem \ref{LGS-13n2-1}.
\begin{Corollary}\label{LKTS-LGS1}
Let $n\equiv 9\pmod{12}$.
\begin{itemize}
\item [$(1)$] If there exists an $\mathrm{LKTS}(n)$, then there exists a $\mathrm{TOC}_{3}(n,5,3)$;\vspace{-0.2cm}
\item [$(2)$] If there exists an $\mathrm{OLKTS}(n)$, then there exists a $\mathrm{TOC}_{3}(n+1,5,3)$.\vspace{-0.2cm}
\end{itemize}
\end{Corollary}

Combine Theorem \ref{Lkts} with Corollary \ref{LKTS-LGS1} to obtain the following results.

\begin{Theorem}\label{LKTS-LGS2}
 Let $a\equiv 3\pmod{6}$ and $a<100$. Further let $m$ be a nonnegative integer and let $m$ be odd if $a\equiv 3\pmod{12}$ and $m$ be even if $a\equiv 9\pmod{12}$.  Then there exist
\begin{itemize}
\item  a $\mathrm{TOC}_{3}(3^{m}a,5,3)$  if $a\not\in \{57,87,93\}$ and
\item  a $\mathrm{TOC}_{3}(3^{m}a+1,5,3)$ if $a\not\in \{51,69,75\}$.
\end{itemize}
\end{Theorem}

%

\begin{Example}\label{LGS1383}
{\rm  Let $X =[8]\times [3]$ and $\mathcal{G}=\big\{\{(x,{i}):i\in [3]\}:x\in [8]\big\}$. 
We construct an $\mathrm{LGS}(1,3,8,3)$, namely a $\mathrm{TOC}_{4}(8,5,3)$, on the point set $X$ with group set $\mathcal{G}$. The following display $9$ mutually disjoint $\mathrm{GS}(1,3,8,3)$s, where $(x,i)$ is written as $x_{i}$. }
\end{Example}\vspace{-0.2cm}
\begin{center}\vspace{-0.2cm}
\begin{tabular}{l l l l l l l l }
$1_1\text{ }2_1\text{ }3_1$&$1_2\text{ }6_2\text{ }7_2$&$1_3\text{ }5_2\text{ }8_2$&$2_2\text{ }4_3\text{ }8_3$&$2_3\text{ }5_1\text{ }7_3$&$3_2\text{ }4_1\text{ }7_1$&$3_3\text{ }6_3\text{ }8_1$&$4_2\text{ }5_3\text{ }6_1$\\
$1_1\text{ }2_1\text{ }4_1$&$1_2\text{ }3_3\text{ }6_1$&$1_3\text{ }5_2\text{ }7_3$&$2_2\text{ }5_1\text{ }6_3$&$2_3\text{ }3_2\text{ }8_2$&$3_1\text{ }4_2\text{ }7_1$&$4_3\text{ }5_3\text{ }8_1$&$6_2\text{ }7_2\text{ }8_3$\\
$1_1\text{ }2_1\text{ }5_1$&$1_2\text{ }3_3\text{ }8_2$&$1_3\text{ }4_2\text{ }6_2$&$2_2\text{ }6_1\text{ }7_2$&$2_3\text{ }3_1\text{ }4_3$&$3_2\text{ }5_3\text{ }7_1$&$4_1\text{ }7_3\text{ }8_3$&$5_2\text{ }6_3\text{ }8_1$\\
$1_1\text{ }2_1\text{ }6_1$&$1_2\text{ }4_3\text{ }8_3$&$1_3\text{ }5_2\text{ }7_2$&$2_2\text{ }5_3\text{ }8_2$&$2_3\text{ }3_2\text{ }4_1$&$3_1\text{ }5_1\text{ }6_2$&$3_3\text{ }7_1\text{ }8_1$&$4_2\text{ }6_3\text{ }7_3$\\
$1_1\text{ }2_1\text{ }7_1$&$1_2\text{ }4_1\text{ }6_3$&$1_3\text{ }3_1\text{ }5_1$&$2_2\text{ }4_2\text{ }8_3$&$2_3\text{ }5_3\text{ }6_1$&$3_2\text{ }4_3\text{ }7_3$&$3_3\text{ }6_2\text{ }8_1$&$5_2\text{ }7_2\text{ }8_2$\\
$1_1\text{ }2_1\text{ }3_2$&$1_2\text{ }5_2\text{ }8_2$&$1_3\text{ }4_2\text{ }7_2$&$2_2\text{ }6_2\text{ }7_1$&$2_3\text{ }4_1\text{ }8_1$&$3_1\text{ }5_3\text{ }7_3$&$3_3\text{ }6_3\text{ }8_3$&$4_3\text{ }5_1\text{ }6_1$\\
$1_1\text{ }2_1\text{ }4_2$&$1_2\text{ }6_2\text{ }8_1$&$1_3\text{ }3_3\text{ }7_1$&$2_2\text{ }6_3\text{ }7_2$&$2_3\text{ }3_1\text{ }5_2$&$3_2\text{ }4_3\text{ }8_2$&$4_1\text{ }5_3\text{ }6_1$&$5_1\text{ }7_3\text{ }8_3$\\
$1_1\text{ }2_1\text{ }4_3$&$1_2\text{ }3_1\text{ }6_1$&$1_3\text{ }5_2\text{ }7_1$&$2_2\text{ }7_2\text{ }8_2$&$2_3\text{ }3_2\text{ }5_3$&$3_3\text{ }4_2\text{ }8_1$&$4_1\text{ }6_2\text{ }7_3$&$5_1\text{ }6_3\text{ }8_3$\\
$1_1\text{ }2_1\text{ }5_2$&$1_2\text{ }4_1\text{ }7_2$&$1_3\text{ }3_3\text{ }8_1$&$2_2\text{ }4_3\text{ }6_2$&$2_3\text{ }3_2\text{ }7_3$&$3_1\text{ }5_1\text{ }6_3$&$4_2\text{ }5_3\text{ }8_3$&$6_1\text{ }7_1\text{ }8_2$\\
\end{tabular}
\end{center}
Let $\alpha$ be the permutation on $X$ defined by $\alpha(x_{i})= x_{i+1}$, where $x\in [8]$, $i,i+1\in [3]$ with reduction by modulo $3$.
Let $\beta$ be the permutation defined by  $\beta(8_{i})=8_{i}$ for  $i\in [3]$ and $\beta(x_{i})=(x+1)_{i}$ if $x\in [7]$, with $x+1\in [7]$ by modulo $7$.
Under the action of the group generated by $\alpha$ and $\beta$, we can obtain $9\cdot 21=189$ mutually disjoint $\mathrm{GS}(1,3,8,3)$s.
Thus we form an $\mathrm{LGS}(1,3,8,3)$.

\begin{Theorem}\label{LGS_d52}For any odd number $n$, there exist
  a $\mathrm{TOC}_{4}(7^{n}+1,5,3)$ and a  $\mathrm{TOC}_{3}(8^{n}+1,5,3)$.
\end{Theorem}
\proof For any odd $n$ and for any prime power $m$, there exists an $\mathrm{RS}(3,m+1,m^{n} + 1)$ by Theorem \ref{KTS}.
There exists an  $\mathrm{LGS}(1,3,8,3)$ by Example \ref{LGS1383}. Then there exists an $\mathrm{LGS}(1,3,7^{n}+1,3)$ by Theorem \ref{LGS-reslovable2}.
Similarly, an $\mathrm{LGS}(1,3,9,2)$ exists by Theorem \ref {LKTS-LGS2} and hence  an $\mathrm{LGS}(1,3,8^{n}+1,2)$ exists  by Theorem \ref{LGS-reslovable2}.
Then we obtain the desired tilings.
\qed

%

\section{Summary and open problems}

Constant-weight codes  have been the subject of extensive research in coding theory.
In this paper, we focussed our perspective on the entire metric space ${\cal H}_q(n,w)$. This area has received relatively little attention.
A  metric space is a set with structure determined by a notion of distance.
Prominent examples of such distance functions that have received significant attention include the Hamming metric, Lee metric, and Johnson metric.
We considered the metric space $\mathcal{H}_{q}(n,w)$ consisting of all words of length $n$ with weight $w$ over an alphabet $\mathbb{Z}_{q}$, equipped with the Hamming metric.
We researched into  tilings $\mathrm{TOC}_{q}(n,d,w)$ of  $\mathcal{H}_{q}(n,w)$ with mutually disjoint optimal $(n,d,w)_{q}$-codes. The main contributions are summarized as follows.

\begin{itemize}
\item We developed constructions for $\mathrm{TOC}_{q}(n,d,w)$s via  $t$-resolvable Steiner systems, see Theorems \ref{LGS-reslovable2}, \ref{LGS-reslovable4};
 \vspace{-0.2cm}
\item  We presented constructions to produce $\mathrm{TOC}_{q}(n,2w-1,w)$s via almost-regular edge-colorings of $K_{n}^{w}$ with certain properties, see Theorems \ref{LGS-13n2-1}, \ref{LGS-13n2-11};\vspace{-0.2cm}
\item The existence problems of $\mathrm{TOC}_{q}(n,d,w)$s were completely resolved in the following situations: (1) $d=2$, (2) $d=2w$, (3) $w=3,q=2$, (4) $w=3,n=4,d=4$, see Theorems \ref{TOC_n_2_w}, \ref{TOC_n_2w_w_1}, \ref{Boundq2--T}, \ref{n43_n=4};\vspace{-0.2cm}
\item  For specific alphabet size $q\ge 3$,  we obtained many infinite families of $\mathrm{TOC}_{q}(n,d,3)$s for distances $d=3,4,5$, see Theorems   \ref{LGS23ng1}, \ref{n43_n=4},  \ref{d4gprime}, \ref{LKTS-LGS2}, \ref{LGS_d52}.
\vspace{-0.2cm}
\end{itemize}

In the process of investigation, we utilized many types of combinatorial objects with resolvable or partitionable properties. 
Some of the existence problems of these  systems have been  remarkably resolved, for example, large sets of Steiner triple systems and almost-regular edge-colorings of $K_n^w$. As a contrast, we also have some widely open problems, such as  large sets of  Kirkman triple systems and $t$-resolvable Steiner systems even for $t=2$. Although we used stronger configurations to treat tilings sometimes, the tiling problems are difficult in general because of the distance constraints. For example, a large set of generalized Steiner triple systems LGS$(2,3,n,2)$ exist for $n=4$ but cannot exist for  $n=3$, which raises the threshold of possible recursive constructions. Likewise, it is not easy to study the large sets of generalized maximum H-packings  LGMHP$(2,3,n,g)$. However, it is a relief if we consider the tiling problems related with $t=1$ cases, so we propose the following problems as the end of the paper.

\begin{itemize}
\item  Establish the existence of an LGS$(1,w,n,g)$ if $w\mid gn$, which amounts to a 1-factorization of the complete $w$-uniform $n$-partite hypergraph $K_{n\times g}^w$ with distance $2w-1$;\vspace{-0.2cm}
\item  Establish the existence of a $\mathrm{TOC}_{q}(n,2w-1,w)$ if  $w\nmid gn$, which is closely related with an almost-regular edge-coloring of  $K_{n\times g}^w$ with distance $2w-1$;
    \vspace{-0.2cm}
    \item Develop constructions of $g$-good almost-regular edge-colorings of the complete hypergraph $K_{n}^{w}$. For small $g$, an asymptotic existence result is promising;\vspace{-0.2cm}
\item Study the existence of $g^{*}$-good almost-regular edge-coloring of $K_{n}^{w}$;\vspace{-0.2cm}
\item Solve the above problems with focus on $w=3$.\vspace{-0.2cm}
\end{itemize}

\section*{Acknowledgments}
The authors are very grateful to Prof. Tuvi Etzion for many helpful suggestions on this topic.



\clearpage

\end{document}